\documentclass[11pt,a4paper]{amsart}
\usepackage[]{amsmath,amsfonts,amsthm,amssymb}

\usepackage{pdfsync}
\usepackage{multicol}
\usepackage{enumerate}
\usepackage[english]{babel}
\usepackage[all]{xy}
\xyoption{all}
\swapnumbers

\usepackage{amsmath}
\usepackage{amssymb}
\usepackage[bbgreekl]{mathbbol}
\usepackage{array}
\usepackage{bbm}
\usepackage{mathrsfs}
\usepackage{dsfont}
\usepackage{amsbsy}
\input xy
\xyoption{all}

\swapnumbers

\theoremstyle{theorem}
\newtheorem{theorem}{Theorem}[section]
\newtheorem{proposition}[theorem]{Proposition}
\newtheorem{lemma}[theorem]{Lemma}

\newtheorem{example}[theorem]{Example}

\theoremstyle{definition}

\newtheorem{definition}[theorem]{Definition}

\newcommand{\aaa}{\ensuremath{\mathcal{A}}}
\newcommand{\bbb}{\ensuremath{\mathcal{B}}}

\newcommand{\fff}{\ensuremath{\mathcal{F}}}

\newcommand{\kkk}{\ensuremath{\mathcal{K}}}

\newcommand{\ppp}{\ensuremath{\mathcal{P}}}

\newcommand{\ttt}{\ensuremath{\mathcal{T}}}

\newcommand{\LL}{\ensuremath{\mathfrak{L}}}

\newcommand{\uu}{\ensuremath{\mathfrak{U}}}

\renewcommand{\a}{\ensuremath{\alpha}}
\renewcommand{\b}{\ensuremath{\beta}}
\renewcommand{\c}{\ensuremath{\gamma}}
\renewcommand{\d}{\ensuremath{\delta}}

\newcommand{\ve}{\ensuremath{\varepsilon}}

\renewcommand{\i}{\ensuremath{\iota}}

\newcommand{\m}{\ensuremath{\mu}}
\newcommand{\n}{\ensuremath{\nu}}

\renewcommand{\r}{\ensuremath{\rho}}
\newcommand{\s}{\ensuremath{\sigma}}
\renewcommand{\t}{\ensuremath{\tau}}

\def\rw{\rightarrow}

\newcommand{\R}{\ensuremath{\mathbb{R}}}
\renewcommand{\P}{\ensuremath{\mathbb{P}}}
\newcommand{\N}{\ensuremath{\mathbb{N}}}
\newcommand{\Q}{\ensuremath{\mathbb{Q}}}
\newcommand{\E}{\ensuremath{\mathbb{E}}}

\newcommand{\frl}{\mathfrak{l}}
\newcommand{\fru}{\mathfrak{u}}

\newcommand{\cat}[1]{\mathsf{#1}}

\renewcommand{\P}{\cat{P}}

\newcommand{\Top}{\ensuremath{\mathsf{Top}}}

\newcommand{\app}{\ensuremath{\mathsf{App}}}

\newcommand{\qmet}{\ensuremath{q\mathsf{Met}}}

\newcommand{\met}{\ensuremath{\mathsf{Met}}}

\renewcommand{\inf}{\myinf}
\renewcommand{\min}{\mymin}

\DeclareMathOperator{\cl}{cl}

\DeclareMathOperator*{\myinf}{in\vphantom{p}f}

\DeclareMathOperator*{\mymin}{mi\vphantom{p}n}

\renewcommand{\phi}{\varphi}

\def\twoheaddownarrow{\rlap{$\downarrow$}\raise-.5ex\hbox{$\downarrow$}}
\def\twoheaduparrow{\rlap{$\uparrow$}\raise.5ex\hbox{$\uparrow$}}

\renewcommand{\phi}{\varphi}  

\allowdisplaybreaks
\numberwithin{equation}{section}

\usepackage{pdfsync}
\begin{document}

\title{Normality in terms of distances and contractions}
\author{E. Colebunders, M. Sioen, W. Van Den Haute}
\maketitle

\vspace{.5cm}

{ \small \noindent Keywords: Normality; Urysohn contractive map; Kat\v{e}tov-Tong interpolation; \\Tietze's contractive extension;  Approach space; (Quasi-)Metric space; Non-expansive map.\\

\vspace{.5cm}

\begin{abstract}

The main purpose of this paper is to explore normality in terms of distances between points and sets. We prove some important consequences on realvalued contractions, i.e. functions not enlarging the distance, showing that as in the classical context of closures and continuous maps, 
normality in terms of distances based on an appropriate numerical notion of \emph{$\c$-separation of sets}, has far reaching consequences on real valued contractive maps, where the real line is endowed with the Euclidean metric. We show that normality is equivalent to (1) separation of $\c$-separated sets by some \emph{Urysohn contractive map}, 
(2) to \emph{Kat\v{e}tov-Tong's interpolation}, stating that for bounded positive realvalued functions, between an upper and a larger lower regular function, there exists a contractive interpolating map and (3)
to \emph{Tietze's extension theorem} stating that certain contractions defined on a subspace can be contractively extended to the whole space.

The appropriate setting for these investigations is the category of approach spaces, but the results have (quasi)-metric counterparts in terms of non-expansive maps. 
Moreover when restricted to topological spaces, classical normality and its equivalence to separation by a Urysohn continuous map, to Kat\v{e}tov-Tong's interpolation for semicontinuous maps and to Tietze's extension theorem for continuous maps are recovered.

\end{abstract}

\section{Introduction}
Normality of a metric space $(X,d)$ in terms of continuous maps is a well known and very useful concept in analysis. In a metric space $(X,d)$ (1) disjoint closed sets can be separated by a continuous Urysohn function, (2) for realvalued bounded functions $g \leq h$ with $g$ upper semicontinuous and $h$ lower semicontinuous, there exists a continuous map $f$ satisfying $g \leq f \leq h$ and (3) a realvalued continuous map defined on some closed subspace has a continuous extension to the whole space. These statements are equivalent to the normality of the underlying topological space $(X, \ttt_d)$ and are known as separation by Urysohn continuous maps (1), as Kat\v{e}tov-Tong's interpolation (2) and as Tietze's extension theorem (3), see for example \cite{BB},\cite{E} \cite{T}, \cite{K}, \cite{Kc}. 

However for some applications in analysis, like for instance the theory of differential equations or fixed point theory, metric spaces with Lipschitz type functions or non-expansive maps are more natural. Such isometric settings get more and more attention like for instance in the study of  approximation by Lipschitz functions in \cite{GJ}, of cofinal completeness and UC-property in \cite{BG}, in investigations on hyperconvexity in \cite{KKO} and on the non-symmetric analogue of the Urysohn metric space in \cite{KS} and \cite{KS2}. For other applications the larger context of approach spaces with contractions is even more suitable as was recently shown in the context of probability m easures \cite{BLC11}, \cite{BLC13} and \cite{BHSC}, or complexity analysis \cite{CWL} and \cite{CWS}.

In this paper we will work in the category $\app$ of approach spaces with contractions \cite{LOW}. The objects of $\app$ are sets $(X,\d)$ endowed with a numerical \emph{distance} $\d(x,A)$ between sets and points (see \ref{distance} for the exact formulation of the axioms) and a map
$f: (X,\delta_X) \rightarrow (Y,\delta_Y)$ is a \emph{contraction} if 
$\forall x\in X, \  \forall A \subseteq X,   \ \delta_Y (f(x), f(A))\leq \delta_X (x,A)$.  $\app$ contains $\met,$ the category of extended pseudometric spaces with non-expansive maps, as well as its non-symmetric counterpart, the category $\qmet$ of extended quasi-pseudometric spaces with non-expansive maps, as fully embedded subcategories, where for a (quasi) metric $d$ the associated approach structure is the natural distance $\d_d(x,A) = \inf_{a \in A} d(x,a)$. Also the category $\Top$ of topological spaces with continuous maps is fully embedded in $\app$ where the distance $\d_{\ttt}(x,A)$ associated to a topology $\ttt$ takes only two values, depending on whether the point $x$ belongs to the closure of $A$ or not.  Moreover every approach space $(X,\d)$ has an underlying quasi-metric $d_\d$ as well as an underlying topological structure $\ttt_\d$ with closures denoted by $\cl(A) = A ^{(0)}$ for $A \subseteq X,$ see \ref{underlying}.

The main purpose of this paper is to show what the meaning and consequences of normality are in terms of distances and contractions. 
Normality for an approach space $(X,\d)$ \ref{normal} is based on an appropriate numerical notion of \emph{$\c$-separation of sets} and on the concept of \emph{contractive scale}, which will be presented in \ref{sep} and \ref{scale}. \emph{Normality} for an approach space $(X,\d)$ states that for every two $\c$-separated subsets $A,B \subseteq X$ there is a contractive scale ``separating'' the closures $A^{}$ and $B^{}$ in the underlying topology. 
In \ref{Ury} separation by \emph{Urysohn contractive maps} (1) states that for every two sets $A,B$ that are $\c$-separated for some $\c>0,$ there exists a contraction $f$ on $(X,\d)$ to the interval $([0, \c],\d_{d_\E})$ with $d_\E$ the Euclidean metric, taking the value $\c$ on the closure $A^{}$ and $0$ on the closure $B^{}$ and this is shown to be equivalent to normality. Based on this equivalence in \ref{Jones} we prove a counterpart of \emph{Jones's Lemma} on separable normal spaces \cite{E}.

The appropriate Kat\v{e}tov-Tong interpolation by contractions \ref{Kat}, is based on bounded functions in the classes $\LL$ of \emph{lower regular} and $\uu$ of \emph{upper regular} functions for an approach space $(X,\d).$ These are the classes of contractions to $[0,\infty],$ endowed with the quasi-metric $d_\P(x,y) = (x - y) \vee 0$ and its dual $d_\P^-$ respectively
 \cite{LOW}, which we recall in \eqref{lower},\eqref{upper}. \emph{Kat\v{e}tov-Tong's interpolation} (2) states that for bounded functions to $[0,\infty]$ satisfying $g \leq h$ with $g$ upper regular and $h$ lower regular, there exists a contractive map $f: (X,\d) \rw ([0,\infty], \d_{d_\E})$ satisfying $g \leq f \leq h,$ and is equivalent to normality.

Next we deal with the problem of extending contractions. Tietze's extension theorem for contractions is a result for each $\c \in \R^+$ and depends on the notion of \emph{development} of a contraction \cite{LOW} which we recall in 5.1. For an approach space $(X,\d)$, a subspace $Y \subseteq X$ and $\c \in \R^+$ we first determine a specific subclass of contractions $f$ on $Y$ to $([0,\infty], \d_{d_\E})$ having a development satisfying a certain condition \ref{Tietze}. \emph{Tietze's extension theorem} (3) states that for all maps $f$ in the particular subclass, there exists a contraction $g$ on $X$ to $([0,\infty], \d_{d_\E})$ of which the restriction to $Y$ coincides with $f$.

The main results in this paper deal with the equivalence of normality with each of the conditions in (1), (2) or (3). The theorems have (quasi)-metric counterparts in terms of non-expansive maps. 
Moreover when restricted to topological (approach) spaces classical normality and its equivalence to the classical separation by a Urysohn continuous map, Kat\v{e}tov-Tong's interpolation for semicontinuous maps and Tietze's extension theorem for continuous maps are recovered.
In section 6 we compare normality to approach frame normality of $\LL$ as studied in \cite{VO} and in section 7 we investigate preservation of normality by maps and subspaces. Although in general compact Hausdorff approach spaces need not be normal, in section 9 we show that the \v{C}ech Stone compactification of $(\b \N, \d_{d_\E})$ is normal.

Categorical considerations linking our normality notions to $\bbbeta \P_+$-normality and to $\bbbeta \P_+$-regularity, as introduced in the context of Monoidal Topology \cite{MT} is work in progress that will be published elsewhere \cite{CSH}.

\section{Preliminaries on approach spaces}

For more details on concepts and results on approach spaces we refer to \cite{LOW} or \cite{apox}.  We recall terminology and basic results that will be needed in this paper.

Usually an extended quasi-pseudometric on a set $X$ is a function $q : X \times X \rightarrow [0, \infty]$ which vanishes on the diagonal and satisfies the triangular inequality and if $q$ vanishes on the diagonal and satisfies both the triangular inequality and symmetry then it is called an extended pseudometric.
So in this paper all such $q : X \times X \rightarrow [0, \infty]$ are allowed to take the value $\infty$ and both distances between two different points can be zero. \emph{From now on, for simplicity in terminology we drop the words ``extended'' and ``pseudo''}, so in this respect our terminology differs from what is commonly used. It is however conform with the terminology in \cite{BLC11}, \cite{CWS} and \cite{CWL}, \cite {LOW}, \cite{MT}. We denote by $\cat{qMet}$ the category of all \emph{quasi-metric spaces with non-expansive maps} as morphisms and by $\cat{Met}$ the full subcategory of all \emph{metric spaces}.

A \emph{distance} on a set $X$ is a function
\begin{equation}\label{distance}
\delta:X\times 2^X\to [0, \infty]
\end{equation} 
with the following properties:
\begin{enumerate}
\setlength{\itemindent}{-3pt}
\item[(D1)] $\forall x\in X,\ \delta (x,\{x\})=0$,
\item[(D2)] $\forall x\in X,\ \delta (x,\emptyset)=\infty$,
\item[(D3)] $\forall x\in X, \ \forall A,B\in 2^X, \ \delta (x,A\cup B)=\min\{\delta (x,A),\delta (x,B)\}$,
\item[(D4)] $\forall x\in X, \ \forall A\in 2^X, \ \forall\ve\in [0,\infty], \ \delta (x,A)\leq\delta (x,A^{(\ve)})+\ve$\\
 with  the enlargement $$A^{(\ve)}=\{x|\delta (x,A)\leq\ve\}.$$
\end{enumerate}
A pair $(X, \d)$ consisting of a set $X$ endowed with a distance $\d$  is called an \emph{approach space}.
For $A \subseteq X$ we denote $\d_A: X \rw [0,\infty]$ the function defined by $\d_A(x) = \d(x,A).$ For $\omega < \infty$ the truncated function $\d_A \wedge \omega$ is denoted by $\d_A^\omega.$

Morphisms between approach spaces are called contractions. A map $f: (X,\delta_X) \rightarrow (Y,\delta_Y)$ is a \emph{contraction} if 
\begin{equation} \label{contract}
\forall x\in X, \  \forall A \subseteq X,   \ \delta_Y (f(x), f(A))\leq \delta_X (x,A)
\end{equation}
The category of approach spaces and contractions
is denoted by ${\mathsf{App}}.$
Contractivity can also be characterised by 
\begin{equation} \label{level}
f(A^{(\ve)}) \subseteq (f(A))^{(\ve)}
\end{equation}
for all $A\subseteq X$ and $\ve\in \R^+,$ where on the lefthandside the enlargement of $A$ depends on $\d_X$, and on the righthandside the enlargement of $f(A)$ on $\d_Y.$

Given an approach space $(X,\d)$ we will also sometimes use the \emph{core operator} $(\i^\omega)_{\omega < \infty}$ instead of the distance
\begin{equation} \label{core}
\i^\omega: X \times 2^X \rw [0,\infty] : (x,A) \mapsto \omega \ominus \d(x,X \setminus A)
\end{equation}
and define the function $\i_A^\omega: X \rw [0,\infty]$ by $\i_A^\omega(x) = \i^\omega(x,A).$ Here on $[0,\infty]$ we shortly write $$x \ominus y:= (x - y) \vee 0.$$ 
For more details on the core we refer to \cite{SH}.

The category ${\mathsf{App}}$ constitutes a framework wherein other important categories can be fully embedded.
${\mathsf{Top}}$ is embedded as a full concretely reflective and concretely coreflective subcategory and $\cat{qMet}$ is embedded as a concretely coreflective subcategory.
The embedding of topological spaces is determined by associating with every topological space $(X, \mathcal{T})$ (with closure of $A$ written as
 $\textup{cl}A$) the distance
$$
\delta_\ttt(x, A) = \begin{cases} 0 & x \in \textup{cl}A \\ \infty & x \not \in \textup{cl}A.\end{cases}
$$
The property of being $\{0, \infty\}$-valued actually characterises approach spaces which are derived from topological spaces.

Every approach space $(X,\delta)$ has two natural topological spaces associated with it, the topological coreflection, which we will also call the \emph{underlying topology}, and the topological reflection. In this paper we will mainly deal with the coreflection which is the topological space $(X,\ttt_{\delta})$ determined by the closure
\begin{equation}\label{underlying}
x \in \textup{cl}A\Leftrightarrow \delta(x, A) =0  \Leftrightarrow x \in A^{(0)}.
\end{equation}

The embedding of quasi-metric spaces is given in the usual way that one defines a distance between points and sets in a metric space, given $q$ we put
$
\delta_q(x, A) = \inf_{a \in A} q(x,a)
.$
We consider two quasi-metrics on $[0,\infty]$, the quasi-metric $d_\P(x,y) = x \ominus y$ and its dual $d_\P^-$ and note that for the Euclidean metric we have $d_\E = d_\P \vee d_\P^{-}$.

An important object in $\cat{App}$ is the space $\mathbb{P},$ as it is an initially dense object in $\app.$ 
Let $\mathbb{P} = ([0,\infty], \delta_{\mathbb{P}})$ where for $x\in [0,\infty]$ and $A \subseteq [0,\infty],$
\begin{equation}\label{defP}
\delta_{\mathbb{P}}(x, A) = \begin{cases} x \ominus \textup{sup} A  & A \neq \emptyset \\ \infty &A = \emptyset. \end{cases}
\end{equation}
On $[0,\infty]$ the structure $\delta_{\mathbb{P}}$ is neither generated by a topological, nor by a quasi-metric structure.
However when $\delta_{\mathbb{P}}$ is induced on $[0,\infty[$ it is generated by the quasi-metric $d_\P.$

The following concepts will play an important role in the sequel.
For an approach space $(X,\d)$ the classes $\LL$ of \emph{lower regular} and $\uu$ of \emph{upper regular} functions, are defined by
\begin{equation}\label{lower} 
\LL = \{f: (X,\d) \rw ([0,\infty], \delta_{\mathbb{P}})\ | \ \text{contractive} \},
\end{equation}
and
\begin{equation}\label{upper}  
 \uu = \{f: (X,\d) \rw ([0,\infty], \d_{d_\P^-})\ | \ \text{bounded, contractive}  \}.
 \end{equation}
Both are stable for taking finite suprema and infima, $\LL$ moreover is stable for arbitrary suprema, and $\uu$ is stable for arbitrary infima.
We denote by $\LL_b$ the bounded functions in $\LL$ and in view of the remark made after \eqref{defP}
we have
\begin{equation}\label{lowerd} 
\LL_b = \{f: (X,\d) \rw ([0,\infty], \d_{d_\P })\ | \ \text{bounded, contractive} \}.
\end{equation}

The collection of all contractions $f: (X,\d) \rw ([0,\infty], \d_{d_\E})$ is denoted by $\kkk$ or more explicitely $\kkk((X,\d))$ or $\kkk((X,\d),([0,\infty], \d_{d_\E}))$ and of all bounded contractions by $\kkk_b.$ Sometimes the structures on domain or codomain will not be mentioned explicitely.
Then we have \cite{SH} 
\begin{equation}\label{intersect}
f \in \uu \cap \LL_b \Leftrightarrow f\in \kkk_b.
\end{equation}

If $\LL$ is the lower regular function frame then the function $\frl: [0,\infty]^X \rw [0,\infty]^X $ defined by
$$
\frl(\mu) := \bigvee\{ \nu \in \LL | \nu \leq \mu\}
$$ is called the \emph{lower hull operator}. This operator is idempotent, monotone, preserves finite infima and for a constant function $\a$ we have $\frl(\mu + \a) = \frl(\mu) + \a.$
The distance can be recovered from the lower hull operator by
\begin{equation}\label{delta}
\d_A = \frl(\theta_A)
\end{equation}
for $A \subseteq X,$ where we use the notation  
$$
\theta_A: X \rw [0,\infty]: x \mapsto \begin{cases} 0  & x \in A \\ \infty & x \not \in A  \end{cases}
$$
and $\theta^\omega_A $ for the truncated function $\theta_A \wedge \omega.$ From \eqref{delta} it follows that $\d_A$ is lower regular.

If $\uu$ is the upper regular function frame then the function $\fru: [0,\infty]^X_b \rw [0,\infty]^X_b $ defined by
$$
\fru(\mu) := \bigwedge\{ \nu \in \uu | \mu \leq \nu\}
$$ is called the \emph{upper hull operator}. This operator is idempotent, monotone, preserves finite suprema and for a constant function $\a$ we have $\fru(\mu + \a) = \fru(\mu) + \a.$
The core can be recovered from the upper hull operator by
\begin{equation}\label{iota}
\i_A^\omega = \fru(\theta_A^\omega)
\end{equation}
for $A\subseteq X, \omega<\infty.$ In particular $\i_A^\omega$ is upper regular.

Approach spaces can be isomorphically described by regular functions or by hulloperators \cite{LOW}, but this will not be needed in this paper. 

\section{Normality and separation by Urysohn maps}
In this section normality is a statement in terms of distances, based on the notions $\c$-separation of sets and contractive scale for two $\c$-separated sets. We prove that contractive scales are in correspondence with contractions and this observation is the cornerstone of our main theorem in section 3 on separation of $\c$-separated sets by Urysohn contractive maps.
\subsection{$\c$-separation of sets and contractive scales}
The basic concept we will use in the definition of normality is $\c$-separation of two sets. It replaces the topological notion that sets have disjoint closures.

\begin{definition}\label{sep}
Let $(X,\d)$ an approach space and $\c>0$. Two sets $A,B \subseteq X$ are called \emph{$\c$-separated} if $A^{(\a)} \cap B^{(\b)}  = \emptyset$, whenever $\a \geq 0$, $\b \geq 0$ and $\a +\b < \c$. 
\end{definition}

A very useful result in section 4, is the impact of $\c$-separation of sets on the regular functions  $\d_A^\c$ and $\i_{X \setminus B}^\c$.

\begin{proposition}\label{inequal}
Let $(X,\d)$ be an approach space. The following are equivalent for subsets $A$ and $B$ and $\c > 0$.
\begin{enumerate}
\item $A$ and $B$ are $\c$-separated 
\item $\i^\c_{X \setminus B} =\c - \d^{\c}_B \leq \d^{\c}_A$
\item $\i^\c_{X \setminus A} =\c - \d^{\c}_A \leq \d^{\c}_B.$
\end{enumerate}
\begin{proof}
(1) $\Rightarrow$ (2): Suppose $A$ and $B$ are $\c$-separated  and let $x \in X$. Either $\d^{\c}_A(x) \geq \c$. In this case we have $\c - \d^{\c}_B(x) \leq \c  \leq \d^{\c}_A(x).$ Or $\d^{\c}_A(x) < \c.$ So for $\d(x,A) = \a$ we have $\a < \c$ and then in view of (1) for $\ve >0$ we have $x \not \in B^{(\c - \a -\ve)}$. So $\c - \d_B^\c(x)< \a + \ve$ and by arbitrariness of $\ve$ we can conclude that $\c - \d_B^\c(x)\leq \a = \d^{\c}_A(x).$\\
(2) $\Rightarrow$ (1): Let $\a, \b \geq 0, \a +\b < \c$ and assume that $x \in A^{(\a)}$. Then we have $\d_A(x) \leq \a < \c$. By (2) it follows that 
$\c - \d^{\c}_B(x) \leq \d^{\c}_A(x) \leq \a < \c - \b$ so $x \not \in B^{(\b)}$.\\
\end{proof}
\end{proposition}

In order to define normality in $\app$, given $\c$-separated sets, we need a counterpart for the topological situation where $\cl(A) \cap \cl(B) = \emptyset$ implies that an open set $G$ can be found such that $\cl(A) \subseteq G \subseteq \cl(G) \subseteq X \setminus \cl(B)$ and that this process can be repeated with open sets defined for every rational number.

\begin{definition}\label{scale}
Let $(X,\d)$ be an approach space. Let $F: \Q\to 2^X$ such that $\bigcup_{q\in \Q} F(q)=X, \bigcap_{q\in \Q} F(q)=\emptyset.$ Then $F$ is a \emph{contractive scale} if it satisfies $$\forall r,s\in \Q: r<s\Rightarrow  F(r) \  \text {and} \  (X \setminus F(s)) \ \text{are}\  (s-r) \text{-separated}$$
\end{definition}

\begin{definition}\label{normal}
An approach space $(X,\d)$ is said to be \emph{normal} if for all
$ A,B\subseteq X,$ for all $\c>0$ with $A$ and $B$ $\c$-separated, a contractive scale $F$ exists such that
\begin{itemize}
\item[(i)] $\forall q\in \Q_0^-: F(q)=\emptyset;$
\item[(ii)] $A^{}\subseteq \bigcap_{q\in \Q^+} F(q);$
\item[(iii)] $B^{}\cap \bigcup_{r\in \Q^+\cap ]0,\c]} F(r)=\emptyset.$
\end{itemize}

\end{definition}

\subsection{Urysohn contractive maps}
The basic result we need in order to link normality to separation of sets by means of contractive maps, describes a correspondence between contractive scales and contractive maps.

\begin{proposition}\label{scale-contract}
Let $F: \Q\to 2^X$ be a contractive scale on an approach space $(X,\d)$. Then
$$f: (X,\d) \to (\R, \d_{d_\E}): x\mapsto\inf\{q\in \Q\mid x\in F(q)\}$$ is a contraction.

\begin{proof}
First notice that since $F$ is monotone increasing, $f$ is well-defined. We will prove contractivity of $f$ by \eqref{level}.
Let $A\subseteq X$ and $x \in A^{(\ve)}$ for some $\varepsilon\in \R^+$ and let $\ve < \r$ be arbitrary.

Pick $q,r \in \Q$ such that $f(x) - \r< q < f(x) - \ve < f(x) + \ve < r < f(x) + \r,$ then clearly we can choose $\c \in \Q_0 $ with $\ve < \c < \min\{ f(x) - q, r - f(x)\},$ which implies 
$
f(x) \in ]q+\c, r-\c[.
$
From the definition of $f$ we can deduce that 
\begin{equation}\label{incl}
f(F(r) \setminus F(q)) \subseteq [q,r].
\end{equation}
By the assumptions made on $F$ we have that $F(r-\c)$ and $X \setminus F(r)$ are $\c$-separated. 
Since $f(x) < r-\c$ we have $x \in F(r - \c)$ and therefore $$x \not \in (X \setminus F(r))^{(\ve)}.$$
We also have that $F(q)$ and $X \setminus F(q +\c)$ are $\c$-separated. 
Since $f(x) > q+\c$ we have $x \in X \setminus F(q+\c),$ which implies that $x \not \in(F(q))^{(\ve)}.$

So by axiom (D3) of the distance we have 
$\d(x, (X \setminus F(r) ) \cup F(q)) > \ve.$
Since on the other hand $x\in A^{(\ve)}$ and therefore $\d(x,A) \leq \ve,$ we have that there exists a point
\begin{equation}\label{point}
a \in (F(r) \setminus F(q)) \cap A.
\end{equation}
Combining \eqref{incl} and \eqref{point} we have $f(a) \in [q,r] \subseteq ]f(x) - \r, f(x) + \r[,$ so for the Euclidean metric we have $d_\E(f(x), f(a)) < \r$. We can conclude that $\d_{d_\E}(f(x), f(A)) < \r$ and by the arbitrariness of $\r$  we finally have $f(x) \in (f(A))^{(\ve)}.$
\end{proof}
\end{proposition}

We call $f$ the \emph{contraction associated to the contractive scale} $F$.

\begin{proposition}\label{contract-scale}
Let $(X,\d)$ be an approach space and $f:(X,\d) \to (\R, \delta_{d_\mathbb{E}})$ a contraction. Then a contractive scale $F:\Q\to 2^X$ exists such that $f$ is the contraction associated to the contractive scale: For all $x\in X$
$$f(x)=\inf\{q\in \Q\mid x\in F(q)\}.$$

\begin{proof}
Assume $f:(X,\d) \to (\R, \delta_{d_\mathbb{E}})$ is a contraction and define
$$F:\Q\to 2^X: r\mapsto \{f\leq r\}.$$
Then clearly $$f(x)=\inf\{q\in \Q\mid x\in F(q)\}$$ for all $x\in X.$ It is obvious that
$$\bigcup_{q\in \Q} F(q)=X, \bigcap_{q\in \Q} F(q)=\emptyset.$$
Now take $r,s\in\Q$ with $r<s$ and $\alpha,\beta\in\R^+$ with $\alpha+\beta<s-r.$ We prove that $F(r)^{(\alpha)}\cap (X \setminus F(s))^{(\beta)}=\emptyset.$ 

Assume that on the contrary $x\in F(r)^{(\alpha)}\cap (X \setminus F(s))^{(\beta)}$ for some $x\in X.$
Since $f$ is a contraction, the following hold:
$$\delta_{d_\E}(f(x),]-\infty,r])\leq \delta_{d_\E}(f(x),f(F(r)))\leq \delta(x,F(r))\leq \alpha,$$
which implies $f(x)\leq r+\alpha$ and
$$\delta_{d_\E}(f(x),]s,\infty[)\leq \delta_{d_\E}(f(x),f(X \setminus F(s)))\leq \delta(x,X \setminus F(s))\leq \beta$$ 
which implies $f(x)\geq s-\beta.$
This would imply
$s-r\leq \alpha+\beta,$ which is a contradiction.

\end{proof}
\end{proposition}

\begin{definition}\label{Ury}
An approach space $(X,\d)$ satisfies \emph{separation by Urysohn contractive maps} if for every $A,B \subseteq X$ and for every $\c >0$, whenever $A$ and $B$ are $\c$-separated, there exists $f \in \kkk((X,\d), ([0,\c],\d_{d_\E}))$ satisfying $f(a) = \c$ for $a\in A^{}$ and $f(b) = 0$ for $b\in B^{}$. 
\end{definition}
Observe that a contraction $f \in \kkk((X,\d), ([0,\c],\d_{d_\E}))$ is continuous with respect to the underlying topologies, so when it satisfying $f(a) = \c$ for $a\in A^{}$ and $f(b) = 0$ for $b\in B^{}$ then we also have $f(a) = \c$ for $a\in A^{(0)}$ and $f(b) = 0$ for $b\in B^{(0)}$.

\begin{theorem}\label{Ury-norm}
An approach space $(X,\d)$ is normal if and only if it satisfies separation by Urysohn contractive maps.

\begin{proof}
First assume that in $(X,\d)$ two $\c$-separated sets can be separated by a Urysohn contraction. Let $A,B\subseteq X$ be $\c$-separared for some $\c >0.$ 
The  given Urysohn map is a contraction $f$ such that $f(a)=0,  f(b)=\c$ for $a\in A^{}, b\in B^{}$. Then by \ref{contract-scale} $$F:\Q\to 2^X:q\mapsto \{f\leq q\}$$ is a contractive scale with the required properties.

Next assume that $(X,\d)$ is normal and let $A,B\subseteq X$ be $\c$-separared for some $\c >0.$
A contractive scale satisfying (i)-(iii) of \ref{normal} exists and let $f$ be its associated contraction. Then $f\wedge \c$ maps all elements of $A^{}$ to $0$ and of $B^{}$ to $\c$.
\end{proof}
\end{theorem}

\subsection{Jones's lemma}
Jones's lemma is an important tool when dealing with normality in concrete topological examples. The lemma asserts that if a separable normal topological space contains a discrete closed set $L$ of cardinality $|L|$, then $2^{|L|} \leq  2^{\aleph_0}.$
In section 7, when we deal with questions on normality of concrete examples of approach spaces,  we need a similar result in terms of distances. In the next proposition by \emph{separability} of an approach space $(X,\d)$ we mean separability of the underlying topology.

\begin{proposition}\label{Jones}

If a separable normal approach space $(X,\d)$ contains a subset $L\subseteq X$ for which some $\c>0$ exists such that for all $A\subseteq L$ the sets $A$ and $L \setminus A$ are $\c$-separated in $(X,\d)$, then $2^{|L|}\leq 2^{\aleph_0}.$

\begin{proof}
By normality and \ref{Ury-norm}, for every $A\in 2^L$, a function $$f_{A}\in \mathcal{K}((X,\d),([0,\c],\d_{d_\E})) \ \text{exists with} \  f_{A}(A^{})\subseteq\{0\}, f_{A}((L \setminus A)^{})\subseteq\{\c\}.$$ Since by $\c$-separation, the sets $A$ and $L \setminus A$ are closed in the underlying topology of $L,$
$$2^L\to \mathcal{K}((X,\d),([0,\c],\d_{d_\E})): A\mapsto f_{A}$$ is an injection, so we have that
$$2^{|L|}\leq |\mathcal{K}((X,\d),([0,\c],\d_{d_\E}))|.$$
Now let $D$ be a countable subset of $X$ that is dense in $(X,\d)$ for the underlying topology. The map
$$\mathcal{K}((X,\d),([0,\c],\d_{d_\E}))\to \mathcal{K}((D,\d),([0,\c],\d_{d_\E})): f\mapsto f\vert_D$$ is injective since $D$ is dense in $(X,\ttt_\d)$ and $([0,\c],\ttt_{d_\E})$ is Hausdorff. Hence $$2^{|L|} \leq |\mathcal{K}((X,\d),([0,\c],\d_{d_\E}))|\leq 2^{\aleph_0}.$$
\end{proof}
\end{proposition}

\subsection{Topological approach spaces}

We prove that the restriction of normality in $\app$ to the full subcategory $\Top$ gives classical normality.

\begin{proposition}\label{normtop}
Let $(X,\ttt)$ be a topological space. The following properties are equivalent
\begin{enumerate}
\item $(X,\ttt)$ is normal in the topological sense
\item$(X,\d_{\ttt})$ is normal in the sense of \ref{normal}.
\end{enumerate}

\begin{proof}
Let $(X,\ttt)$ be a topological space $A,B \subseteq X$ and $\c>0$. By the fact that the associated distance  is two-valued we have 
$$
A,B \ \text{are}\  \c \ \text{-separated for} \  (X,\d_\ttt) \Leftrightarrow \cl(A) \cap \cl(B) = \emptyset.
$$
Moreover by the coreflectivity of $\Top$ in $\app$ we have 
$$f: (X, \d_\ttt) \rw (\R, \d_{d_\E}) \ \text{ is a contraction} \Leftrightarrow f: (X, \ttt) \rw (\R, \ttt_{d_\E})\ \text{is continuous}.$$
\end{proof}
\end{proposition}

Remark that the topological Jones's Lemma can be deduced from \ref{Jones}. For a topological space $(X,\ttt)$ a subset $L$ being closed and discrete, implies that for $\c>0$ arbitrary and for all $A\subseteq L$ the sets $A$ and $L \setminus A$ are $\c$-separated.

\subsection{Examples}
In this section we present some examples of normal approach spaces, among them are some normal quasi-metric approach spaces. Examples of non-normal quasi-metric approach spaces will be encountered in \ref{exInorm}, \ref{exVO} and \ref{Xn}.
\begin{example}\label{P}
The approach space $\mathbb{P} = ([0,\infty], \delta_{\mathbb{P}}))$ introduced in \eqref{defP} and the quasi-metric approach spaces $([0,\infty], \delta_{d_\P})$ and  $([0,\infty], \delta_{d_\P^-})$ are normal.
\begin{proof}
It is easy to see that in  $\mathbb{P}$ two non-empty subsets $A$ and $B$ always have $0 \in A^{(0)} \cap B^{(0)}$. So for $\c>0$ arbitrary, there are no non-empty $\c$-separated sets. Hence normality of $\mathbb{P}$ is trivially fulfilled.

The same argument holds to show normality of  $([0,\infty], \delta_{d_\P})$.  In order to show normality of $([0,\infty], \delta_{d_\P^-})$ we can use $\infty$ instead of $0$.
\end{proof}
\end{example}

Next we look at a quasi-metric space inducing the Sorgenfrey line and show that it is a normal approach space.
\begin{example}\label{q}
The quasi-metric approach space $([0, \infty[,\d_q)$ defined by $q(x,y) =  y - x$ if $x\leq y$ and $q(x,y) = \infty$ if $x >y$ is normal.

\begin{proof}
Let $A,B\subseteq X,$ $\c>0$, and assume that $A$ and $B$ are $\c$-separated for $\d_q.$
First we prove that this implies that $A$ and $B$ are also $\c$-separated for $\d_{d_\E}.$
The enlargements for $q$ will be denoted by a superscript $q$, and for the Euclidean metric by a superscript $\E$.

Assume that $x\in A^{(\m)_\E}\cap B^{(\n)_\E}$ for some $\m+\n<\c.$ Choose $\varepsilon>0$ such that $\m+\n+2\varepsilon<\c.$ Then $a_0\in A, b_0\in B$ exist such that
$$d_\E(x,a_0)<\m+\varepsilon,\quad d_\E(x,b_0)<\n+\varepsilon$$
and hence $d_\E(a_0,b_0)<\c.$
We consider three cases. If $a_0=b_0,$ then clearly $a_0\in A^{(0)_q}\cap B^{(0)_q}$. If $a_0>b_0$, we can change the roles of $A$ and $B$ and continue with the next case.
So assume that $a_0<b_0.$ Set $\alpha=0, \beta=\m+\n+2\varepsilon.$ Then $\alpha+\beta<\c$, $a_0\in A^{(0)_q}$ and
$$\delta_q(a_0,B)=\inf_{b\in B}q(a_0,b)\leq d(a_0,b_0)<\m+\n+2\varepsilon$$
and therefore $$a_0\in A^{(\alpha)_q}\cap B^{(\b)_q}.$$

Since $A$ and $B$ are $\c$-separated for the Euclidean metric we can apply \ref{metric} to find a contraction
$$f\in \mathcal{K}(([0,\infty[,\delta_{d_\E}),([0,\c],\delta_{d_\E}))$$ with $f(A^{})\subseteq\{0\}$ and $f(B^{}\subseteq\{\c\}$. Since $\delta_\E\leq \delta_q,$ we have that $$f\in \mathcal{K}(([0,\infty[,\delta_q),([0,\c],\delta_{d_\E})).$$  
\end{proof}
\end{example}

\section{Kat\v{e}tov-Tong's interpolation}

By a deep and beautiful result in $\Top$ normality can be characterised by means of interpolation between semicontinuous functions. This result is known as Kat\v{e}tov-Tong's result \cite{K}, \cite{T}. In this section we solve the question on what type of interpolation by means of a contractive function should be used in order to catch normality.

\subsection{Interpolation}

\begin{definition}\label{Kat}
An approach space $(X,\d)$ satisfies \emph{Kat\v{e}tov-Tong's interpolation} if 
for bounded functions to $[0,\infty]$ satisfying $g \leq h$ with $g$ upper regular and $h$ lower regular, there exists a contractive map $f: (X,\d) \rw ([0,\infty], \d_{d_\E})$ satisfying $g \leq f \leq h.$
\end{definition}

In \cite{T} Tong proved a general lemma on a lattice $M$ and a sublattice $K.$  The lemma provides sufficient conditions for elements $s\in K_\delta, t\in K_\sigma$ with $s\leq t$  to have an interpolating $u\in K_\sigma\cap K_\delta$ satisfying $s\leq u\leq t,$ where
$K_\sigma=\{\bigvee_{n} t_n\mid \forall n: t_n\in K\}$ and  $K_\delta=\{\bigwedge_{n} t_n\mid \forall n: t_n\in K\}.$ 

In Theorem \ref{Ury-Kat} we will apply Tong's lemma to the special situation $K=\mathcal{K}(X,[0,\omega])$ and $M=[0,\omega]^X$.  In this particular case we have $K_\s \subseteq \LL$ and $K_\d \subseteq \uu$ and the lemma takes the following simpler form.

\begin{lemma}\label{Tong}
Let $K=\mathcal{K}((X,\d),([0,\omega],\d_{d_\E}))$ and $M=[0,\omega]^X,$ let $s\in K_\delta = \{\bigwedge_{n\geq 1} t_n\mid \forall n: t_n\in K\}$ and $t\in K_\sigma = \{\bigvee_{n} t_n\mid \forall n: t_n\in K\}$ with $s\leq t$ then a $u \in K_\sigma\cap K_\delta$ exists satisfying $s\leq u\leq t.$
\end{lemma}

\subsection{Equivalence with normality}
The main theorem of this section links Kat\v{e}tov-Tong's interpolation to separation by Urysohn contractive maps.

\begin{theorem}\label{Ury-Kat}
For an approach space $(X,\d)$, the following are equivalent.
\begin{enumerate}
\item $(X,\d)$ satisfies Kat\v{e}tov-Tong's interpolation.
\item $\forall A,B\subseteq X, \forall\omega<\infty: (\iota^\omega_A\leq \delta^\omega_B \Rightarrow \exists f\in \mathcal{K}_b((X,\d)): \iota^\omega_A \leq f\leq  \delta^\omega_B).$
\item $(X,\d)$ satisfies separation by Urysohn contractive maps
\end{enumerate}
\begin{proof}
That (1) implies (2) follows from \ref{delta} and \ref{iota}.\\
$(2)\Rightarrow(3.)$. Let $A,B\subseteq X, \c>0$ such that $A$ and $B$ are $\c$-separated. Put $\omega=\c$ in (2). Hence, by Proposition \ref{inequal}, $\iota^\c_{A^c}\leq \delta^\c_B.$
By $(2)$ a contraction $f$ exists such that $\iota^\c_{A^c}\leq f\leq \delta^\c_B.$ Clearly $f\in \mathcal{K}(X,[0,\c])$ and $$\forall a\in A^{}: f(a)=\c, \ \forall b\in B^{}: f(b)=0.$$

$(3)\Rightarrow(1)$. Let $\varphi\in \uu, \psi\in \LL$ with $\varphi\leq  \psi\leq \omega$ for some $\omega\in \R.$ For $k,m,n\in \N$ with $m\leq k<n$, set
$$A_{m,n}=\{\psi\leq \omega \frac{m}{n}\}, \quad B_{k,n}=\{\varphi \geq \omega (\frac{k}{n}+\frac{1}{2n})\}.$$
We show that $A^{(\a)}_{m,n}\cap B^{(\b)}_{k,n} = \emptyset$ for all $\a + \b < \c $ with $\c= \omega \frac{2k - 2m +1} {2n}.$

Let $\alpha+\beta< \c.$ If $x\in A_{m,n}^{(\alpha)}$ then since $\psi: X \rw ([0,\omega], d_\P)$ is contractive, by \ref{level} we have $\psi(x) \in [0,wm/n]^{(\a)}$. It follows that $\inf_{z \leq \omega m/n} (\psi(x) \ominus z) \leq \a$ and therefore 
$$\psi(x) \leq \frac{\omega m}{n} + \a.$$
Similarly, assuming $x \in B^{(\b)}_{k,n}$ and using the contractivity of $\varphi$ to the codomain endowed with $d^-_\P$  it follows that 
$$\varphi(x) \geq \frac{\omega (2k +1)}{2n} - \b.$$
Since $\varphi\leq \psi,$ the assumtion $x\in A_{m,n}^{(\alpha)}\cap B^{(\b)}_{k,n}$ would imply $\frac{\omega (2k +1)}{2n} - \b \leq \frac{\omega m}{n} + \a$ which is impossible.

By (3), for $1<m\leq k<n$,
a contraction $$f_{m,n}^k \in \kkk(X, [\omega \frac{m+1}{n}, \omega \frac{2k+3}{2n} \wedge \omega])$$ exists with $f_{m,n}^k\vert_{A_{m,n}}=\omega \frac{m+1}{n}$ and $ f_{m,n}^k\vert_{B_{k,n}}=\omega \frac{2k+3}{2n} \wedge \omega.$
Define the contractions $$f_{m,n}=\bigvee_{k=m}^{n-1}f_{m,n}^k \ \text{and} \  f_n = \bigwedge _{m=2}^{n-1} f_{m,n}.$$

Next we show that $\varphi \leq f_n$ whenever $n \geq 3$.
Let $x\in X$ and $1<m<n$, either $x \not\in B_{m,n}$, then 
$$\varphi(x) \leq \omega \frac{2m +1}{2n} \leq \omega \frac{m+1}{n} \leq f_{m,n}(x),$$ 
or $x\in B_{m,n}$, then we again consider two cases. If $x \in B_{n-1,n}$ then
$$\varphi(x)\leq \omega = f_{m,n}^{(n-1)}(x) \leq f_{m,n}(x).$$
Otherwise, a minimal $k$ exists with $m< k \leq n-1$ and $x\in B_{k-1,n}, x \not \in B_{k,n}.$ Then we have
$$\varphi(x) < \omega \frac{2k+1}{2n} = \omega \frac{2(k-1)+3}{2n} = f_{m,n}^{(k-1)}(x) \leq f_{m,n}(x).$$

Next we show that $\bigwedge_{n\geq 3} f_n \leq \psi$. In order to do so,
 for $x\in X$ we prove that $$f_n(x) - \psi(x) \leq \omega \frac{2}{n}$$ for every $n \geq 3$. 
 Fix $x\in X$, then one of three possibilities holds. First if $f_n(x)\leq \psi(x),$ we are done. Secondly if $x\notin A_{m,n}$ for all $m \geq 2$, then $\psi(x)>\omega(n-1)/n$. Since $f_n(x)\leq \omega,$ we have that $$f_n(x)-\psi(x)\leq \omega\frac{1}{n}<\omega\frac{2}{n}.$$
Thirdly, if some minimal $m\geq 2$ exists such that $x\in A_{m,n},$ then $\psi(x)\geq\omega(m-1)/n$ and $f_{m,n}(x)=\omega(m+1)/n.$ So $$f_{n}(x)-\psi(x)\leq f_{m,n}(x)-\psi(x)\leq\omega \frac{2}{n}.$$
So we can conclude that
$$
\varphi \leq \inf_{n\geq 3} f_n \leq \psi \ \text{and therefore also} \ \omega - \psi \leq \omega - \inf_{n\geq 3} f_n  \leq \omega.
$$ 
Remark that in the last inequality we have  $\omega - \psi \in \uu$ and $\omega - \inf_{n\geq 3} f_n \in \LL.$ So we can repeat our argument above with $\varphi$ replaced by $\omega - \psi $ and $\psi$ replaced by $\omega - \inf_{n\geq 3} f_n.$ We find contractions $g_n$ for $n \geq 3$ satisfying $\omega - \psi \leq \inf_{n \geq 3} g_n \leq \omega - \inf_{n\geq 3} f_n.$ Thus we have
$$\varphi\leq \bigwedge_{n\geq 3} f_n\leq \omega-\bigwedge_{n\geq 3} g_n=\bigvee_{n\geq 3} (\omega-g_n)\leq \psi.$$
Using Lemma \ref{Tong} an $f\in L_\sigma\cap L_\delta\subseteq \LL_b\cap \uu = \kkk_b$ exists with $\bigwedge_n f_n\leq f \leq \bigvee_n (\omega-g_n).$ Hence $f$ is a contraction with $$\varphi\leq f\leq \psi.$$
\end{proof}
\end{theorem}

\subsection{Topological and metric approach spaces}First we show that when theorem \ref{Ury-Kat} is restricted to topological spaces, we recover the classical theorem on interpolation between semicontinuous maps.

\begin{proposition}\label{Kattop}
Let $(X,\ttt)$ be a topological space. The following properties are equivalent
\begin{enumerate}
\item $(X,\ttt)$ satisfies the classical topological Kat\v{e}tov-Tong's interpolation theorem for bounded (positive) realvalued maps.
\item $(X,\d_\ttt)$ satisfies Kat\v{e}tov-Tong's interpolation in the approach sense.
\end{enumerate}

\begin{proof}
Let $(X,\ttt)$ be a topological space, then in view of the coreflectivity of $\Top$ in $\app$ a map $f: (X,\d_\ttt) \rw ([0,\infty], \d_{d_\P})$ is contractive if and only if $f: (X,\ttt) \rw ([0,\infty], \ttt_{d_\P})$ is continuous. The topology $ \ttt_{d_\P}$ is generated by $\{ ]x,\infty] | x\geq 0\},$ so $\LL_b$ coincides with the class of bounded lower semicontinuous maps. In the same way $\uu$ coincides with the class of bounded upper semicontinuous maps. Moreover the class of bounded contractions $f: (X,\d_\ttt) \rw ([0,\infty], \d_{d_\E})$ coincides with the class of bounded continuous maps to $([0,\infty], \ttt_\E),$ with $\ttt_\E$ the Euclidean topology.
So in the topological case, Kat\v{e}tov-Tong's interpolation in the approach sence, is equivalent to some topological Kat\v{e}tov-Tong's interpolation theorem for bounded maps.
\end{proof}
\end{proposition}
Remark that when \ref{normtop} and \ref{Kattop} are combined, we can conclude that the restriction in (1) of \ref{Kattop} to bounded maps is not really needed.

Next we investigate the normality properties restricted to the full subcategory $\qmet$. For a  quasi-metric space $(X,d)$ a function $f: (X,\d_d) \rw ([0,\infty], \d_{d_\E})$ is a contraction if and only if $f: (X,d) \rw ([0,\infty], {d_\E})$ is non-expansive, so both Urysohn separation and Kat\v{e}tov-Tong interpolation deal with non-expansive maps in that context. 

Next we focus on metric spaces.
We have the following result which makes use of the fact proved \cite{LOW} that on a metric approach space, by the symmetry,  each of the classes $\LL_b$ and $\uu$ coincide with $\kkk_b$.

\begin{proposition}\label{metric}
For an approach space $(X,\d_d)$ with $(X,d)$ a metric space we have
\begin{enumerate}
\item If $A,B\subseteq X$ are $\c$-separared for some $\c >0$ then $f = (\d_d)_B \wedge \c$ is a Urysohn non-expansive map taking values $0$ on $B^{(0)}$ and $\c$ on $A^{(0)}.$
\item If $g \leq h$ with $g$ upper regular and $h$ lower regular, then $f = g$ (or $f=h$) is a non-expansive map inbetween, so Kat\v{e}tov-Tong's interpolation is trivially fulfilled
\end{enumerate}
and therefore $(X,\d_d)$ is normal.

\begin{proof}
(1) As we recalled in \eqref{delta} the map $(\d_d)_B: (X, d) \rw [0,\infty]: x \mapsto \d_d(x,B)$ is lower regular. Since  in the metric case $\LL_b$ and $\uu$ coincide, we can conclude that $(\d_d)_B: (X,d) \rw ([0,\infty], d_\E)$ is non-expansive. Clearly it takes value $0$ on $B^{(0)}$ and since $A,B$ are $\c$-separated, $(\d_d)_B \geq \c$ on $A^{(0)}$.\\ 
(2) Immediately from $\LL_b =\uu =\kkk_b$.

\end{proof}
\end{proposition}

That the reverse implication does not hold follows from example \ref{q}, where we describe a quasi-metric approach space that is not a metric one, but nevertheless is normal.

\section{Tietze's extension theorem}
The topological version of Tietze's theorem ensures continuous extensions of continuous maps defined on closed subsets. In order to prove an approach version that is equivalent to normality, it is not possible to ensure contractive extensions for all contractive maps defined on subsets that are closed in the underlying topology. The resulting property would be too strong.

In the proof of the main theorem \ref{Tietze-Ury} of this section, given a normal approach space an extension of a given contraction, satisfying certain conditions, is constructed by first defining an upper regular extension and a lower regular extension and then applying Kat\v{e}tov-Tong's theorem \ref{Ury-Kat}. This is done firstly for functions taking only a finite number of values and then for contractions, by  so called developments. This means that for a bounded function some uniform approximation from below by functions taking only a finite number of values is needed. This approximation technique was described in \cite{LOW} and is briefly recalled below.

\subsection{Extensions based on developments}

Given a bounded $f \in [0,\infty]^X$ a family $(\mu_\ve)_{\ve>0}$ of functions taking only a finite number of values, written as $$\left(\mu_\varepsilon:=\bigwedge_{i=1}^{n(\varepsilon)} \left(m_i^\varepsilon + \theta_{M_i^\varepsilon}\right)\right)_{\ve>0} \text{with} \ (M^\ve_i)_{i=1}^{n(\ve)} \ \text{a partitioning of}\  X$$  and all $m_i^\ve \in \R^+,$ for $\ve >0,$ is called a \emph{development of $f$} if for all $\ve >0$
$$
\mu_\ve \leq f \leq \mu_\ve + \ve.
$$
In 1.1.30 of \cite{LOW} it is shown that in this case the lower regular hull of $f$ can be calculated as 
\begin{equation} \label{lower}
\frl(f)=\bigvee_{\ve>0} \frl(\mu_\ve) = \bigvee_{\ve>0} \bigwedge_{i=1}^{n(\varepsilon)} \left(m_i^\varepsilon + \frl (\theta_{M_i^\varepsilon})\right).
\end{equation} 
From 1.1.32 and 1.1.33 \cite{LOW} it also follows that for the upper regular hull of $f$ we have 
\begin{equation}\label{upper}
\fru(f) = \bigvee_{\ve>0} \fru(\mu_\ve).
\end{equation} 

If $(X,\d)$ is an approach space and $Y \subseteq X$ then the induced distance on $Y$ is the restriction of $\d: X \times \ppp(X) \rw [0,\infty]$ to $Y \times \ppp(Y)$. We will still denote it by $\d$.

\begin{definition}\label{Tietze}
An approach space $(X,\d)$ satisfies \emph{Tietze's extension theorem } if for every $Y \subseteq X$ and $\c\in \R^+,$ and every $f \in \mathcal{K}((Y,\d),([0,\c], 
\d_{d_\E}))$, having a development $\left(\mu_\varepsilon:=\bigwedge_{i=1}^{n(\varepsilon)} \left(m_i^\varepsilon + \theta_{M_i^\varepsilon}\right)\right)_{0<\varepsilon<1}$ such that
$$\forall x\notin Y, \forall\varepsilon\in ]0,1[, \forall 1\leq l,k\leq n(\varepsilon):m_l^\varepsilon-m_k^\varepsilon\leq \delta_{M_k^\varepsilon}(x)+ \delta_{M_l^\varepsilon}(x),$$
there exists a contractive extension
$$ g\in \mathcal{K}((X,\d),([0,\c], \d_{d_\E})): g\vert_Y=f.$$

\end{definition}

\begin{theorem}\label{Tietze-Ury}
The following conditions on an approach space $(X,\d)$ are equivalent
\begin{enumerate}
\item Normality
\item Tietze's extension theorem.
\end{enumerate}

\begin{proof}
(1) $\Rightarrow$ (2):  Let $Y\not = \emptyset$ and $\c$ be fixed and suppose $f \in \mathcal{K}(Y,[0,\c])$ has a development $\left(\mu_\varepsilon:=\bigwedge_{i=1}^{n(\varepsilon)} \left(m_i^\varepsilon + \theta_{M_i^\varepsilon}\right)\right)_{0<\varepsilon<1}$ such that
$$\forall x\notin Y, \forall\varepsilon\in ]0,1[, \forall 1\leq l,k\leq n(\varepsilon):m_l^\varepsilon-m_k^\varepsilon\leq \delta_{M_k^\varepsilon}(x)+ \delta_{M_l^\varepsilon}(x).$$ Without loss of generality we may assume that all $M_i^\ve$ are nonempty subsets of $Y$.
For $\varepsilon \in ]0,1[,$ we have that
$$\mu_\varepsilon:=\bigwedge_{i=1}^{n(\varepsilon)} \left(m_i^\varepsilon + \theta_{M_i^\varepsilon}\right)\leq f\leq \mu_\varepsilon+\varepsilon.$$
Set $\omega=\c+1.$ Then by pointwise verification we see that on $Y$ for $\varepsilon\in ]0,1[$
$$\mu_\varepsilon=\bigvee_{i=1}^{n(\varepsilon)} \left(\theta^\omega_{(Y \setminus M_i^\varepsilon)} \ominus (\omega- m_i^\varepsilon)\right).$$ 
By \eqref{upper}  and \eqref{iota} for the upper regular hull of $f$ on $Y$ we obtain 
\begin{eqnarray*}
\fru_Y(f) =  \bigvee_{0<\varepsilon<1} \fru(\mu_\ve)
 &=& \bigvee_{0<\varepsilon<1}\bigvee_{i=1}^{n(\varepsilon)} \left(\fru(\theta^\omega_{(Y\setminus M_i^\varepsilon)} \ominus (\omega- m_i^\varepsilon)\right) \\
&=& \bigvee_{0<\varepsilon<1}\bigvee_{i=1}^{n(\varepsilon)} \left(\iota^\omega_{(Y \setminus M_i^\varepsilon)} \ominus (\omega- m_i^\varepsilon)\right)
\end{eqnarray*}
Using the fact that $f$ is contractive on $Y$ we know that the lower en upper regular hull $\frl_Y(f)$ and $\fru_Y(f)$ coincide, moreover by \eqref{lower} and \eqref{delta}
$$\frl_Y(f)=\bigvee_{0<\varepsilon<1} \bigwedge_{i=1}^{n(\varepsilon)} \left(m_i^\varepsilon + \delta_{M_i^\varepsilon}\right) 
= \bigvee_{0<\varepsilon<1}\bigvee_{i=1}^{n(\varepsilon)} \left(\iota^\omega_{(Y \setminus M_i^\varepsilon)} \ominus (\omega- m_i^\varepsilon)\right)=\fru_Y(f).$$

Next we define extensions on $X$.  Let 
$$\hat{\mu}: X \rw [0,\infty] \ \text{ be the function with} \ \hat{\mu}(x) = f(x),  x \in Y \ \text{and} \  \hat{\mu}(x) = \c, x \not \in Y,$$ 
and $$\check{\mu}: X \rw [0,\infty] \ \text{ be the function with} \  \check{\mu}(x) =f(x), x \in Y\ \text{and} \  \check{\mu}(x) = 0, x \not \in Y.$$ 

Similarly for all $0<\ve<1$ let 
$$\hat{\mu}_\ve: X \rw [0,\infty] \ \text{with} \ \hat{\mu}_\ve(x) = \mu_\ve(x), x \in Y \ \text{and} \  \hat{\mu}_\ve(x) = \c, x \not \in Y$$
 and 
 $$\check{\mu}_\ve: X \rw [0,\infty] \ \text{with} \ \check{\mu}_\ve(x) = \mu_\ve(x), x \in Y \ \text{and} \  \check{\mu}_\ve(x) = 0, x \not \in Y.$$ 
Applying the regular hulls on $X$ we obtain 
$$\frl(\hat{\mu_\varepsilon})=  \frl(\bigwedge_{i=1}^{n(\ve)} (m_i^\ve + \theta_{M_i^\ve}) \wedge (\c + \theta_{(X \setminus Y)} ))= \bigwedge_{i=1}^{n(\varepsilon)} \left(m_i^\varepsilon + \delta_{M_i^\varepsilon}\right) \wedge (\c + \d_{(X \setminus Y)})$$
and
$$\fru(\check{\mu_\varepsilon})= \fru \left(\bigvee_{i=1}^{n(\ve)}( \theta^\omega_{(X \setminus M_i^\ve)} \ominus (\omega - m_i^\ve))\right) = \bigvee_{i=1}^{n(\varepsilon)} \left(\iota^\omega_{(X \setminus M_i^\varepsilon)} \ominus (\omega- m_i^\varepsilon)\right).$$
In view of $$\hat{\mu}_\ve \leq \hat{\mu} \leq \hat{\mu}_\ve + \ve \ \text{and} \  \check{\mu}_\ve \leq \check{\mu} \leq \check{\mu}_\ve + \ve$$ for the regular hulls of $\hat{\mu}$ and $\check{\mu}$ on $X$ we get 
$$ \frl(\hat{\mu})=\bigvee_{0< \ve <1} \frl(\hat{\mu}_\ve) = \left(\bigvee_{0 < \ve <1} \bigwedge_{i=1}^{n(\varepsilon)} (m_i^\varepsilon + \delta_{M_i^\varepsilon})\right) \wedge (\c + \d_{(X \setminus Y)})$$ 
and
$$ \ \fru(\check{\mu})=\bigvee_{0<\ve<1} \fru(\check{\mu}_\ve) = \bigvee_{0<\ve<1} \bigvee_{i=1}^{n(\varepsilon)} \left(\iota^\omega_{(X \setminus M_i^\varepsilon)} \ominus (\omega- m_i^\varepsilon)\right).$$ 

We now verify the claim that $\fru(\check{\mu}) \leq \frl(\hat{\mu})$ on $X$.
By pointwise verification, in case $x \in Y$ we easily obtain $$f(x) = \check{\mu}(x)  \leq \fru(\check{\mu})(x) \leq f(x) = \frl(\hat{\mu})(x)$$ and hence  
\begin{equation}\label{onY}
\fru(\check{\mu})(x) = f(x) = \frl(\hat{\mu})(x). 
\end{equation}
In case $x \in X \setminus Y$ for $\varepsilon\in ]0,1[$, we have that
$$\frl(\hat{\mu_\varepsilon})(x)=\bigwedge_{i=1}^{n(\varepsilon)} \left(m_i^\varepsilon + \delta_{M_i^\varepsilon}\right) \wedge \c$$
and
$$\fru(\check{\mu_\varepsilon})(x)=\bigvee_{i=1}^{n(\varepsilon)} \left(\iota^\omega_{(X \setminus M_i^\varepsilon)} \ominus (\omega- m_i^\varepsilon)\right).$$
For $1\leq i\leq n(\varepsilon),$ we also have that
$$\iota^\omega_{(X \setminus M_i^\varepsilon)} \ominus (\omega- m_i^\varepsilon)=(\omega\ominus \delta_{M_i^\varepsilon})\ominus (\omega- m_i^\varepsilon)=m_i^\varepsilon\ominus \delta_{M_i^\varepsilon}.$$
Now choose $1\leq k,l \leq n(\varepsilon)$ such that $m_k^\varepsilon + \delta_{M_k^\varepsilon}(x)$ is minimal and \\$\iota^\omega_{(X \setminus M_l^\varepsilon)}(x) \ominus (\omega- m_l^\varepsilon)$ is maximal. Then
$$m_k^\varepsilon + \delta_{M_k^\varepsilon}(x)\geq \iota^\omega_{(X \setminus M_l^\varepsilon)}(x)-(\omega- m_l^\varepsilon)=m_l^\varepsilon- \delta_{M_l^\varepsilon}(x)$$ by assumption. So the claim follows.

By Kat\v{e}tov-Tong's extension, there exists $g \in \kkk_b(X)$ such that 
$$0 \leq \fru(\check{\mu}) \leq g \leq \frl(\hat{\mu}) \leq \hat{\mu} \leq \c$$
which takes values in $[0,\c]$ and by \eqref{onY} satisfies $g\vert_Y=f.$

(2) $\Rightarrow$ (1): Let $A,B\subseteq X,$ $\c > 0$ and assume that  $A$ and $B$ are $\c$-separated. Set $Y=A^{(0)}\cup B^{(0)}$ and define
$$f:Y=Y^{(0)}\to [0,\c]: x\mapsto\begin{cases}
\c & x\in A^{(0)}\\
0 & x\in B^{(0)}.\end{cases}$$
We show that $f$ is a contraction. Let $x\in Y, C\subseteq Y.$ If $x\in A^{(0)}$ and $C\cap A^{(0)}\neq \emptyset$ or $x\in B^{(0)}$ and $C\cap B^{(0)}\neq \emptyset$, then $\delta_\E(f(x),f(C))=0.$ So assume that $C\subseteq A^{(0)}$ and $x\in B^{(0)}.$ Then $$\delta_\E(f(x),f(C))=\delta_\E(\c,\{0\})=\c.$$ Since $A^{(0)} \cap B^{(0)} = \emptyset$, we have that $\delta(x,A^{(0)})\geq \c$ and hence $\delta(x,C)\geq \c.$ 

The function $f$ on $Y$ has a development where $$\mu_\varepsilon=\theta_{A^{(0)}}\wedge (\c+\theta_{B^{(0)}})$$ for all $\varepsilon$. In order to check that this development satisfies the condition in Tietze's extension theorem, fix $x\not \in Y,$ $\ve>0$ and $l,k\in \{1,2\}.$ Then  $$m_l-m_k\leq \delta_{M_k}(x)+ \delta_{M_l}(x)$$ trivially holds for all cases except for $k=1,l=2.$ In this case, $\delta_{A}(x)+\delta_{B}(x) < \c$ would imply $x \in A^{(\a)} \cap B^{(\b)}$ for $\a = \d_{A}(x)$ and $\b = \d_{B}(x)$.  It follows that 
$$\c \leq \delta_{A}(x)+\delta_{B}(x).$$

Applying (2) choose an extension  $g\in \mathcal{K}(X,[0,\c])$ of $f$. Then $g$ is the required Urysohn function.
\end{proof}
\end{theorem}

\subsection{Topological approach spaces}
First we show that when theorem \ref{Tietze-Ury} is restricted to topological spaces, we recover the classical theorem on extensions of continuous maps.

\begin{proposition}
Let $(X,\ttt)$ be a topological space then the following are equivalent
\begin{enumerate}
\item $(X,\d_\ttt)$ has Tietze's extension theorem for contractions.
\item $(X,\ttt)$ satisfies the classical topological Tietze theorem for continuous maps.
\end{enumerate} 
\begin{proof}
(1) $\Rightarrow$ (2): Let $Y \subseteq X$ be closed and $f: (Y, \ttt_Y) \rw ([0,1], \ttt_\E)$ a continuous map. Then by coreflectivity of $\Top$ in $\app$ we have $f: (Y, \d_{\ttt_Y}) \rw ([0,1], \d_{d_\E})$ is contractive. For $x \not \in Y$ we have $x \not \in \cl_X(M)$ and hence $\d_{\ttt_X} (x,M) = \infty$, for every subset $M$ of $Y.$ It follows that every development of $f$ satisfies the condition in \ref{Tietze}. Hence there exists a contraction $g: (X, \d_{\ttt}) \rw ([0,1], \d_{d_\E})$ with $g | Y = f.$ Clearly $g: (X, \ttt) \rw ([0,1], \ttt_\E)$ is continuous.\\
(2) $\Rightarrow$ (1): By (2) the topological space $(X,\ttt)$ is normal in the topological sense. By \ref{normtop} $(X, \d_\ttt)$ is normal in the approach sense and by \ref{Tietze} Tietze's extension theorem is fulfilled.
\end{proof}
\end{proposition}

\section{Linking to other types of normality}

In this section we investigate the relation of normality in $\app$ to other existing notions of normality. Normality has been studied for frames \cite{PP} and based on this definition ``approach frame normality'' was investigated in \cite{VO} for approach frames. Given an approach space $(X,\d)$ the lower regular function frame $\LL$ is an approach frame and when applied to this particular approach frame one obtains a condition of approach frame normality of $\LL,$ of which an equivalent expression is formulated in (3) of \ref{Inorm}. We prove that this condition is srictly weaker than normality introduced in \ref{normal}.

We also compare normality of $(X,\d)$ with normality of the underlying topology and prove that in fact both notions are unrelated.

In the context of Monoidal Topology \cite{MT} normality and regularity for lax algebras are introduced and studied.  Approach spaces have been described as lax algebras, namely as objects of $\bbbeta \P_+$-$\mathsf{Cat}$ \cite{CH} and so $\bbbeta \P_+$-normality and $\bbbeta \P_+$-regularity, as defined in \cite{MT}, apply to $\app$.
Categorical considerations linking our normality notions to $\bbbeta \P_+$-normality and to $\bbbeta \P_+$-regularity is work in progress that will be published elsewhere \cite{CSH}.

\subsection{Approach frame normality}

For the study of approach frame normality we refer to \cite{VO} or \cite{LO}, where the equivalent characterisation we are using under (3) was proved.
\begin{proposition}\label{Inorm}
Let $(X,\d)$ be an approach space. Consider the following properties
\begin{enumerate}
\item $(X,\d)$ is normal
\item For $A, B \subseteq X,$ $\c$-separated for some $\c>0,$ there exists $C \subseteq X$ such that 
$A$ and $C$ are $\c/2$-separated and $X \setminus C$ and $B$ are $\c/2$-separated.
\item $\LL$ is approach frame normal: For $A, B \subseteq X,$ $\ve>0 $ such that\\ $A^{(\ve)} \cap B^{(\ve)} = \emptyset $ there exist $\r>0,$ $C \subseteq X$ with 
$$
A^{(\r)} \cap C^{(\r)} = \emptyset \ \text{and} \ (X \setminus C)^{(\r)} \cap B^{(\r)} = \emptyset.
$$
\end{enumerate}
Then we have (1) $\Rightarrow$ (2) $\Rightarrow$ (3)

\begin{proof}
(1) $\Rightarrow$ (2): Suppose $(X,\d)$ is normal and $A, B \subseteq X$ are $\c$-separated for some $\c>0.$ By \ref{Ury} we have a Urysohn contraction $f$  with $f(A^{}) \subseteq \{0\}$ and $f(B^{}) \subseteq \{\c\}.$
Let $C = \{ f > \c / 2\}$ and let $\s,\t\geq 0$ with $\s + \t < \c / 2.$ Choose $\s',\t'$ with $\s < \s',$ $\t < \t'$ and $\s' + \t' < \c / 2.$
We claim that $$
A^{(\s)} \cap C^{(\t)} = \emptyset \ \text{and} \ (X \setminus C)^{(\s)} \cap B^{(\t)} = \emptyset.
$$
Otherwise if there is some $x \in A^{(\s)} \cap C^{(\t)}$ then there would exist $a \in A$  and $c \in C$ satisfying $f(x) < \s', f(c) > \c / 2$ and $| f(x) - f(c) | < \t'$, which is incompatible with $\s' + \t' < \c / 2.$ 

Or if there is some $z \in (X \setminus C)^{(\s)} \cap B^{(\t)}$ then there would exist $b \in B$ and $u \in (X \setminus C),$ with $\c \leq f(b), | f(z) - f(b) | < \t'$ and $f(u) \leq  \c / 2, | f(z) - f(u) | < \s'.$ This again is incompatible with $\s' + \t' < \c / 2.$ \\
(2) $\Rightarrow$ (3):  Let $A, B \subseteq X,$ $\ve>0 $ such that $A^{(\ve)} \cap B^{(\ve)} = \emptyset,$  then clearly $A$ and $B$ are $\ve$-separated. So by (2) there exists $C \subseteq X$ such that for all $\s,\t\geq 0$ with $\s,\t\geq 0$ with $\s + \t < \c / 2.$ We have
$$
A^{(\s)} \cap C^{(\t)} = \emptyset \ \text{and} \ (X \setminus C)^{(\s)} \cap B^{(\t)} = \emptyset.
$$
Then $\r = \ve / 5$ fulfills the condition in (3).
\end{proof}
\end{proposition}

None of the implications in \ref{Inorm} is reversible. This is shown by the next examples. 

\begin{example}\label{exInorm}
Consider the following quasi-metric space $(X,d)$ and the associated approach space $(X, \d_d)$.
Let $X= \{x,y,z\}$ with $d(x,z) = 1$ , $d(y,z) =2$ and $d(x,y)= 4$. All other distances are infinite, except for $d(a,a) = 0$ whenever $a \in X$.
We prove the following:
\begin{enumerate}
\item $X$ is not normal
\item $X$ has property (2) in \ref{Inorm}
\end{enumerate}

\begin{proof}
(1) We use the characterization of normality by Urysohn separating contractions. Let $A = \{x\} ,B = \{y\},\c =4.$ For  $\a +\b < \c$ we clearly have $\{x\}^{(\a)} \cap \{y\}^{(\b)} = \{x\} \cap \{y\}= \emptyset.$  

If $f$ would be a contraction $f: (X, \d_d) \rw [0,\c]$, with values $\c$ on $A$ and $0$ on $B,$ then we would have
$|f(x) -f(z)| \leq 1$ which implies $f(z)\geq 3$ and $|f(y) -f(z)| \leq 2$ which implies $f(z) \leq 2$. It is clear that such a contraction does not exist.\\
(2) Suppose $A$, $B$ are $\c$-separated for $\c >0$. Observe that in case $\c \leq 2$ we have $\c /2 \leq 1$ and then it suffices to find $C$ with $A \cap C = \emptyset$ and $(X \setminus C) \cap B = \emptyset.$ 
We have this staightforward situation $\c \leq 2$ in all cases except in the one
where $\{x\}$ and $\{y\}$ are $\c$-separated. Then $\c \leq 4$. Take $C =\{y\}$ and assume $\c +\d < \c / 2 \leq 2$. Then clearly $\{x\}^{(\c)} \cap \{y\}^{(\d)} = \{x\} \cap \{y\} = \emptyset$ and $\{x,z\}^{(\c)} \cap \{y\}^{(\d )}= \{x,z\}^{(\c)} \cap \{y\} = \emptyset.$
\end{proof}
\end{example}

The following example was used in \cite{MT} for other purposes.

\begin{example}\label{exVO}
Let $X=(\{x,y,z,w\},\delta_d)$ with $d(x,z)=d(y,z)=d(w,z)=1, d(w,x)=d(w,y)=2$ and $d(a,a)=0, d(a,b)=\infty$ for all other $a\neq b$. The approach space $(X,\d_d)$ satisfies (3) but not (2) in \ref{Inorm}.
\begin{proof}
For $A, B \subseteq X,$ $\ve>0 $ such that $A^{(\ve)} \cap C^{(\ve)} = \emptyset $ let $C=B$ then  
$$
A \cap C = \emptyset \ \text{and} \ (X \setminus C) \cap B = \emptyset.
$$
Clearly for $\r = 1/2$  the condition in (3) is fulfilled.

Set $A=\{x\}, B=\{y,z,w\}.$ Then clearly they are $4$-separated. Now assume $C$ is some subset of $X.$ Either $z\in C$, but then
$C^{(1)}=X$ so with $\s = 0, \t =1$ the condition in (2) fails, or $z\notin C,$ then with $\s =1, \t = 0$ the condition fails.

\end{proof}
\end{example}

\subsection{Normality of the underlying topology}

Every approach space $(X,\d)$ has an underlying topological space $(X, \ttt_\d).$ In this subsection we comment on possible implications between  normality of the topological space $(X, \ttt_\d)$ and normality of the approach space $(X,\d).$ 
As will become clear both concepts are unrelated.

Example \ref{exInorm} is a space that was shown not to be normal in the approach sense. However its underlying topology is discrete. 
That also the other implication is false, is shown by the next example.
\begin{example}\label{qS}
We start from the quasi-metric approach space in \ref{q} on $[0, \infty[$ defined by $q(x,y) =  y - x$ if $x\leq y$ and $q(x,y) = \infty$ if $x >y$. The quasi-metric approach space $X =([0,\infty[ \times [0,\infty[, q_S)$ where 
$$
q_S(a,b) = q(a',b') + q(a'',b'')
$$
with $a=(a',a'')$ and $b=(b',b'')$ is normal whereas its underlying topology coincides with the Sorgenfrey plane and hence is not normal.

\begin{proof}
The topological coreflection of $X$ is the topological space $$([0,\infty[ \times [0,\infty[, \ttt_q \times \ttt_q)$$ which is known not to be normal.
To prove normality of $X$ by Urysohn separation, we consider the smaller metric
$$
d_S(a,b) = d_\E(a',b') + d_\E(a'',b'')
$$
on $[0,\infty[ \times [0,\infty[$ and 
for $A,B \subseteq [0,\infty[ \times [0,\infty[$ we claim that 
$$A^{(\a)_{q_S}} \cap B^{(\b)_{q_S}} = \emptyset \ \text{whenever}\  \a + \b < \c \Rightarrow A^{} \cap B^{(\eta)_{d_S}} = \emptyset \ \text{whenever} \  \eta < \c.$$
Suppose on the contrary $A^{(\a)_{q_S}} \cap B^{(\b)_{q_S}} = \emptyset$ whenever $\a + \b < \c$  but there exists $\eta < \c$, with $A^{} \cap B^{(\eta)_{d_S}} \not= \emptyset.$ Choose $\sigma$ such that $\eta < \sigma < \c$ and $a \in A^{}$ and $b \in B$ with $d_S(a,b) < \s$. There are different cases for the position of $b$ with respect to $a$.\\
In case $a'\leq b'$ and $a''\leq b''$ we have $q_S(a,b) = d_S(a,b) < \s$ and $a\in A^{} \cap B^{(\s)_{q_S}}$ which is impossible.\\
In case $b'\leq a'$ and $b''\leq a''$ we have $q_S(b,a) = d_S(b,a) < \s$ and $b\in A^{(\s)_{q_S}} \cap B^{(0)_{q_S}}$ which is impossible.\\
In case $b'\leq a'$ and $a''\leq b''$ for $(b',a'')$ we put $q_S((b',a''),a) = a'-b'=\a$ and $q_S((b',a''),b) = b''-a''=\b.$ Then we have $\a + \b = d_\E(a',b') + d_E(a'',b'') =d_S(a,b)< \s < \c$ and $(b',a'') \in A^{(\a)_{q_S}} \cap B^{(\b)_{q_S}}$ which is impossible.\\
In case $a'\leq b'$ and $b'' \leq a''$ for $(a',b'')$ we put $q_S((a',b''),a) = a''-b''=\a$ and $q_S((a',b''),b) = b'-a'=\b.$ Then we have $\a + \b = d_\E(a'',b'') + d_\E(a',b') =d_S(a,b)< \s < \c$ and $(a',b'') \in A^{(\a)_{q_S}} \cap B^{(\b)_{q_S}}$ which is impossible.\\
 This proves our claim.

Assume that $A$ and $B$ are $\c$-separated for $q_S$ and consider the bounded lower regular function $$\d_{d_S}(\cdot, B) \wedge \c : ([0,\infty[ \times [0,\infty[, d_S) \rw ([0,\infty],d_\E).$$ Since the domain is a metric space, the function $\d_{d_S}(\cdot, B) \wedge \c$ is a contraction, which remains contractive when $[0,\infty[ \times [0,\infty[$ is endowed with $q_S$. In view of the claim above for $x \in A^{}$ we have $\d_{d_S}(x,B) > \eta$ for every $\eta < \c$, hence also $\d_{d_S}(x,B) \geq \c.$ So the contraction $\d_{d_S}(\cdot, B) \wedge \c$ has value $\c$ on 
$A^{}.$ On the other hand $\d_{d_S}(\cdot, B) \wedge \c$ has value $0$ on $B^{}.$ 
\end{proof}
\end{example}

\section{Preservation of normality}

\noindent Since $\Top$ is concretely reflective in $\app$, it is clear that the fact that normality is not (finitely) productive in $\Top$ implies that normality is not (finitely) productive in $\app$ either.
However with respect to subspaces and maps in general there are some nice results in $\Top,$ as there
normality is preserved under closed continuous surjections and it
is preserved under taking closed subspaces. We obtain some preservation results in the context of $\app$.

\subsection{Preservation by maps}
We recall the definition of open and closed expansive maps in $\app$  in terms of  regular functions and refer to 1.4.1  and 1.4.3 in \cite{LOW}.

\begin{definition}
A function $f: (X,\d) \rightarrow (X', \d')$ between approach spaces with regular function frames $\LL, \uu$ and $\LL', \uu'$ is said to be
\begin{enumerate}
\item closed expansive if: $f(\m) \in \LL'$ whenever $\m \in \LL,$ with $f(\m)(y) = \inf_{x \in \m^{-1}(y)} \m(x).$
\item open expansive if: $f(\n) \in \uu'$ whenever $\m \in \uu.$
\end{enumerate}
\end{definition}
\begin{proposition}
Normality is preserved by contractive surjections that are at the same time open and closed expansive. 
\begin{proof}
Let $f: (X,\d) \rightarrow (X', \d')$ be a closed and open expansive contractive surjection, between approach spaces with lower regular function frames $\LL$ and $\LL'$. Let $g \leq h$ with $g$ upper regular and $h$ lower regular on $(X',\d')$. Then since f is contractive $g\circ f$ and $h\circ f$ are respectively upper regular and lower regular on $(X,\d)$. So there exists a contraction $\m: (X,\d) \rightarrow ([0,\infty], d_\E)$ satisfying  $g\circ f \leq \m \leq h\circ f.$
Now consider the function $f(\m): (X',\d') \rw ([0,\infty], d_\E)$ which lies between $g$ and $h.$ Moreover since $f$ is closed $f(\m)$ is lower regular and since $f$ is open $f(\m)$ is upper regular. So finally $f(\m)$ is contractive.
\end{proof}
\end{proposition}

\subsection{Preservation by subspaces}
In \cite{LOW} injective closed expansive contractions between approach spaces $f: (X,\d) \rw (X',\d')$ are characterized as embeddings such that $\d'_{f(X)} = \theta_{f(X)}$.

\begin{proposition}\label{emb}
If $f: (X,\d) \rw (X',\d')$ is an injective closed expansive contraction and $(X',\d')$ is normal, then so is $(X,\d)$.
\begin{proof}
Let $A,B$ be $\c$-separated subsets of $(X,\d)$ for some $\c>0.$ Since $$\d'(f(x), f(A)) = \d(x,A)$$ (and similar for $B$) and $\d'(y, f(X)) = \theta_{f(X)}$  we have $(f(A))^{(\a)} \cap (f(B))^{(\b)} = \emptyset,$ whenever $\a + \b < \c$. Let $h'$ be a contraction  $h': (X',\d') \rw ([0,\c], \d_{d_\E})$ satisfying $h'((f(A))^{}) \subseteq \{\c\}$ and $h'((f(B))^{}) \subseteq \{0\}.$ Then $h = h' \circ f$ is a contraction on $(X, \d)$ satisfying $h(A^{}) \subseteq \{\c\}$ and $h(B^{}) \subseteq \{0\}$.
\end{proof}
\end{proposition}

Clearly for an embedding $f: (X,\d) \rw (X',\d')$  in $\app$ the map $f: (X,\ttt_\d) \rw (X', \ttt_{\d'})$ is an embedding in $\Top$ and it is called a closed embedding if moreover $f(X)$ is closed in $(X', \ttt_{\d'})$. Contrary to the result in \ref{emb}, embeddings that are closed in the topological coreflection, do not preserve normality. To build an example of this kind we use the power of Jones's result as established for approach spaces in \ref{Jones}. 

\begin{example}\label{Xn}
We start from the example in \ref{qS}, the normal space $$X =([0,\infty[ \times [0,\infty[, q_S)$$ 
and consider  $n \in \N_0$ fixed. The quasi-metric approach space $X_n$ on $$\{ (x',x'') \in [0,\infty[ \times [0,\infty[ \ |\  x'' \geq -x' +n \},$$ endowed with the trace of $q_S$ is closed in the underlying topology, but is not normal.

\begin{proof}
Clearly $X_n$ has a separable underlying topology since $$X_n \cap \Q \times \Q =\{ (x',x'') \in [0,\infty[ \times [0,\infty[ \cap \Q \times \Q \ |\  x'' \geq -x' +n \}$$ is dense in $X_n.$
Consider the subspace $$L = \{(x',x'') \in [0,\infty[ \times [0,\infty[ \ |\  x'' = -x' +n \}$$  of $X_n$ with the induced quasi-metric $q_S$ and fix $\c = 1$. 

Clearly the restriction of $q_S$ on $L$ takes value $\infty$ on couples of different points, so it is discrete.
Let $(x',x'')\in X_n \setminus L$ and $(a', -a'+n) \in L$ arbitrary different, we either have $x' >a'$ and then $q(x',a') =\infty$ or $x' \leq a'$ and then $x'' > -x' +n \geq -a' +n$ which implies $q(x'', -a' +n) = \infty$. In both cases $q(x',a') + q(x'', -a'+n) = \infty$. 

It follows that for an arbitrary non-empty subset $A \subseteq L,$ the sets $A$ and $L \setminus A$ are $\c$-separated for the
associated distance of $\d_{q_S}.$ 

Assume $X_n$ is normal, then by Jones's result developed for approach spaces in \ref{Jones}, for the cardinality $|L|$ we would have $2^{|L|} \leq 2^{\aleph_0}$. Since $|L| = 2^{\aleph_0}$ a contradiction follows.

\end{proof}
\end{example}

\section{Compact Hausdorff approach spaces}

By compactness of an approach space $(X,\d)$ we mean the underlying topological space $(X, \ttt_\d)$ is compact in the topological sense, and by Hausdorffness we mean the underlying topology is Hausdorff \cite{LOW}.
As example \ref{exInorm} shows, a compact Hausdorff approach space $(X,\d)$ need not be normal (although the underlying topology is normal).

One of the interesting concrete examples developed in the framework of approach theory is the \v{C}ech Stone compactification of $(\N,d_\E).$  The compactification allows to equip $\b \N$ with a canonical approach structure which generates the  \v{C}ech Stone topology as underlying structure and which extends the usual metric $d_\E$ on $\N,$ an impossible situation in the usual setting, since the topological compactification is not metrizable \cite{LOW}, \cite{apox}.

\subsection{The \v{C}ech Stone Compactification of $(\N,\d_{d_\E})$}
In this section we prove that the \v{C}ech Stone compactification $\b^* \N$  in $\app$ is normal.
In order to treat the \v{C}ech Stone compactification $\b^* \N =(\b \N, \d_{\b\N})$ of the approach space $(\N, \d_{d_\E})$ we first have to recall from \cite{LOW} that the topological coreflection of $\b^* \N  = (\b \N, \d_{\b\N})$ is the usual (topological) \v{C}ech Stone compactification of $\N$ with the discrete topology and that the \v{C}ech Stone compactification has the universal property with respect to bounded contractions to $(\R, \d_{d_\E}).$

 As usual $\b \N$ is the underlying set of the topological \v{C}ech Stone compactification of the discrete space $\N$. Points of $\b \N$ are ultrafilters on $\N$ and are denoted by $p,q, \cdots$. For a subset $\aaa \subseteq \b\N$ let $\fff_\aaa = \bigcap_{q \in \aaa} q$. When $\N$ is endowed with the Euclidean metric $d_\E$ the distance of the approach \v{C}ech Stone compactification $\b^* \N$ is given by
$$
\d_{\b\N}(p,\aaa) = \sup_{F \in \fff_\aaa} \inf \{ \a \ |\  F^{(\a)_{d_\E}} \in p \}.
$$
Observe that in view of the fact that  on $\N \times \N$ the function $d_\E$ takes values in $\N$, the function $\d_{\b\N}$ too takes values in $\N$ on $\b \N \times 2^{\b \N}$.

\begin{proposition} \label{beta}
The \v{C}ech Stone compactification $\b^*\N  = (\b \N, \d_{\b\N})$ of the approach space $(\N, \d_{d_\E})$ is normal.

\begin{proof}
For $\aaa$ a subsets of $\b \N$ and for $r$ in $\N$ we denote $\fff_\aaa^{(r)_{d_\E}}$ the filter generated by $\{ F^{(r)_{d_\E}} \ | \ F \in \fff_\aaa \}$. 
Let $\aaa$ and $\bbb$ be subsets of $\b \N$ and $r$ and $s$ in $\N.$
If $p \in \b\N$ satisfies $\fff_\aaa^{(r)_{d_\E}} \subseteq p$ and $\fff_\bbb^{(s)_{d_\E}} \subseteq p$ then $\d_{\b\N}(p,\aaa) \leq r$ and 
$\d_{\b\N}(p,\bbb) \leq s$ and so $p \in \aaa^{(r)_{\b\N}} \cap \bbb^{(s)_{\b\N}}$.

For $\c \in \R^+_0$ assume that $\aaa^{(r)_{\b^*\N}} \cap \bbb^{(s)_{\b^*\N}} = \emptyset$ for all $(r,s) \in \N \times \N$ with $r+s < \c$.
By the previous observation, for each such $(r,s)$ we can choose $F_{(r,s)}\in \fff_\aaa$ and $G_{(r,s)}\in \fff_\bbb$ such that for $\d_{d_\E}$ we have $F^{(r)}_{(r,s)} \cap G^{(s)}_{(r,s)}  = \emptyset$. Put 
$$A = \bigcap_{\{(r,s)|r+s<\c\}} F_{(r,s)} \in \fff_\aaa \ \ \text{and}\ \ B = \bigcap_{\{(r,s)|r+s<\c\}} G_{(r,s)} \in \fff_\bbb.$$ 
Then it follows that for each $(r,s) \in \N \times \N$ with $r+s < \c$ these sets satisfy $A^{(r)} \cap B^{(s)} = \emptyset.$
Since the metric approach space $(\N,\d_{d_\E})$ is normal \ref{metric}, we can find a contraction $$f: (\N,\d_{d_\E}) \rw ([0,\c],\d_{d_\E})$$ satisfying $f(A) \subseteq \{0\}$ and $f(B) \subseteq \{\c\}$. Moreover since $f$ is bounded, by the universal property there exists a unique contractive extension
$$
\hat{f}: (\b \N, \d_{\b\N}) \rw ([0,\c],\d_{d_\E})
$$ 
such that the restriction to $\N$ coincides with $f.$

Since $A \in \fff_\aaa,$ for every $p\in \aaa$ we have $A\in p$ and hence $p \in \cl_{\b\N}A.$ 
 This implies $\cl_{\b\N} \aaa \subseteq \cl_{\b\N} A.$ Using 
 the continuity of $\hat{f}: \b\N  \rw ([0,\c], \ttt_\E)$  we have
$$
\hat{f}(\aaa^{}) \subseteq \hat{f}(\cl_{\b\N}\aaa) \subseteq \hat{f}(\cl_{\b\N}A )\subseteq \cl_{\ttt_\E}f(A) \subseteq \{0\}
$$
and analogously we have $\hat{f}(\bbb^{}) \subseteq \{ \c \}$.

\end{proof}
\end{proposition}

\vspace{.8cm}
\noindent {\small 
E. Colebunders, M. Sioen:
\address{Vakgroep Wiskunde, Vrije Universiteit Brussel, Pleinlaan 2, 1050 Brussel, Belgium}\\
\email{evacoleb@vub.ac.be, msioen@vub.ac.be}\\

\noindent W. Van Den Haute:
\address{Departement Wiskunde-Informatica, Universiteit Antwerpen, Middelheimlaan 1, 2020, Antwerpen, Belgium}\\
\email{wouter.vandenhaute@uantwerpen.be}}

\end{document}